\documentclass{article}

\begin{document}

\title{Volume estimates for K\"ahler-Einstein metrics and rigidity of complex structures}
\author{X-X. Chen and S. K. Donaldson}
\date{\today}
\maketitle


\newcommand{\oZ}{\overline{Z}}
\newcommand{\uf}{\underline{f}}
\newcommand{\tg}{\tilde{g}}
\newcommand{\cA}{{\cal A}}
\newcommand{\bQ}{{\bf Q}}
\newcommand{\cX}{{\cal X}}
\newcommand{\bC}{\mbox{${\bf C}$}}
\newcommand{\bR}{\mbox{${\bf R}$}}
\newcommand{\Var}{{\rm Var}}
\newcommand{\Av}{{\rm Av}}
\newcommand{\Vol}{{\rm Vol}}
\newcommand{\Dim}{{\rm Dim}}
\newcommand{\cO}{{\cal O}}
\newcommand{\cW}{{\cal W}}
\newcommand{\cL}{{\cal L}}
\newcommand{\Tr}{{\rm Tr}}
\newcommand{\Zmax}{Z_{{\rm max}}}
\newcommand{\Zmin}{Z_{{\rm min}}}
\newcommand{\Ch}{{\rm Ch}}
\newcommand{\bP}{\mbox{${\bf P}$}}
\newcommand{\uA}{\mbox{${\underline{A}}$}}
\newcommand{\uM}{\mbox{${\underline{M}}$}}
\newcommand{\um}{\mbox{${\underline{m}}$}}
\newcommand{\ur}{\mbox{${\underline{r}}$}}
\newcommand{\db}{\overline{\partial}}
\newtheorem{Goal}{Goal}
\newtheorem{thm}{Theorem}
\newtheorem{prop}{Proposition}
\newtheorem{lem}{Lemma}
\newtheorem{defn}{Definition}
\newtheorem{cor}{Corollary}
\newcommand{\oK}{\overline{K}}
\newcommand{\dbd}{\sqrt{-1} \partial\overline{\partial}}
\newcommand{\ulambda}{\underline{\lambda}}
\newcommand{\olambda}{\overline{\lambda}}
\newcommand{\Riem}{{\rm Riem}}
\newcommand{\Ric}{{\rm Ric}}

\section{Introduction}

This is a continuation of our previous paper \cite{kn:CD}. Let $M$ be a compact
K\"ahler-Einstein manifold with non-negative scalar curvature and for $r>0$ let $Z_{r}$ be the $r$-neighbourhood of the points where $\vert \Riem\vert \geq r^{-2}$. Our purpose is to estimate the volume of $Z(r)$. In the previous paper we considered manifolds of complex dimension $3$ and here we extend the results to all dimensions (under very slightly different hypotheses). We use a different approach, exploiting the rigidity of complex structures on quotient singularities. This also gives another approach to the three-dimensional case. The basic technique develops results of Tian\cite{kn:T}, for the 3-dimensional case, with the difference that we work with complex domains rather than CR-structures. Meanwhile, as mentioned in \cite{kn:CD}, Cheeger and Naber have posted a preprint \cite{kn:CN} which reaches the same general conclusions using different arguments. In addition, Tian has informed us that he obtained similar results some while ago.

Throughout the paper we will make the following standing assumptions.    

\

\begin{itemize}
\item $(M,g)$ is a compact K\"ahler-Einstein manifold of complex dimension $n$   with $\Ric(g)=\lambda g$, and $\lambda>0$. As usual, we write $\omega$ for the metric 2-form.
\item ${\rm Diam}(M)\leq D $. 
\item The class $[\omega/2\pi]\in H^{2}(M)$ is integral. 
\end{itemize}

Note that these hypotheses have the following standard consequences
\begin{equation} \Vol(M) \geq V_{0}=(4\pi)^{n}/n!. \end{equation}
This just follows from the integrality condition.
\begin{equation} \Vol(B(x,r))\geq \kappa(n,D) r^{2n}; \end{equation}
for metric balls $B(x,r)\subset M$. This follows from Bishop-Gromov comparison.
\begin{equation} {\rm Diam} (M) \geq D_{0}; \end{equation}
we can take $D_{0}$ to be the radius of the Euclidean ball with volume $V_{0}$, again by Bishop comparison. 
\begin{equation}  \lambda \leq \lambda_{0}; \end{equation} 
This follows from the diameter bound and Myers' Theorem.

\

Our main result is a \lq\lq small energy'' estimate. Recall that for a ball $B(x,r)\subset M$ with $r\leq {\rm Diam}(M)$ we define the normalised energy
$$    E(x,r)= r^{2-n} \int_{B(x,r)} \vert \Riem \vert^{2}.
$$   

Then we have

\begin{thm} There are $\epsilon_{0}>0, K$ such that if 
and $E(x,r)\leq \epsilon_{0}$ then $\vert {\rm Riem} \vert \leq K r^{-2}$ on $B(x,r/2)$.
\end{thm}
 Given this we obtain, using just the same line of argument as in \cite{kn:CD}
\begin{cor}
$${\rm Vol}(Z_{r})\leq C(n,D) E(M) r^{4}$$
where $E(M)$ is the square of the $L^{2}$ norm of the curvature.
\end{cor}
Recall here that $E(M)$ is a topological invariant, determined by the Chern classes of $M$ and the K\"ahler class. 

For completeness we also mention that we get an \lq\lq approximate monotonicity'' property for the normalised energy.

\begin{cor} For every $\epsilon>0$ there is a $\delta>0$ so that if $E(x,r)\leq \delta$ then for any  $r'\leq r/2$ and $y\in B(x,r/2)$ we have $E(y,r')\leq \epsilon$. 
\end{cor}

This follows easily from Theorem 1.

\

We will deduce Theorem 1 from the following result.
\begin{thm}
Let $(M_{i}, g_{i})$ satisfy the conditions above; let $x_{i}\in M_{i}$ for each $i$ and let $l_{i}$ be any sequence of numbers which tends to infinity.
Suppose $(M_{i},l_{i}^{2} g_{i}, x_{i})$ has based Gromov-Hausdorff limit
$M_{\infty}$. Then $M_{\infty}$ is not a product $\bC^{n-q}\times \bC^{q}/\Gamma$ where $q>2$ and $\Gamma\subset U(q)$ acts freely on $S^{2q-1}$.

\end{thm}

In turn, Theorem 2 will be proved entirely by complex geometry. The main ingredient is a result on rigidity of complex structures which may have independent interest and, as far as we are aware, is new (in the case when $n>q$). This uses some recent work of Chakrabarti and Shaw  \cite{kn:CS}.
Given $q\leq n$ and a real number  $a>1$, let $V(a)\subset \bC^{n}$ be the domain
$$   \{ z,w)\in {\bf C}^{n-q}\times {\bf C}^{q}: \vert z\vert <a, a^{-1}< \vert w \vert<a\}. $$
Let $\Gamma\subset U(q)$ be as above (acting freely on $S^{2q-1}$) and write
$V_{\Gamma}(a)$ for the quotient of $V(a)$ by $\Gamma$. Fix any $a_{1}, b_{1}$ with  $b_{1}<a_{1}$ so
$V_{\Gamma}(b_{1})\subset V_{\Gamma}(a_{1})$.

\begin{thm} Let $q\geq 3$ and $J'$ be a deformation of the standard complex structure $J$ on $V_{\Gamma}(a_{0})$. If the deformation is sufficiently small in $C^{1,\alpha}$ then there is a diffeomorphism from $V_{\Gamma}(b_{1})$ to a domain in $V_{\Gamma} (a_{1})$ which pulls back $J'$ to the standard complex structure on $V_{\Gamma}(b_{1})$  and which is close in $C^{2,\alpha}$ to the inclusion map.   
\end{thm}

Here of course we mean that the diffeomorphism can be forced as close to the inclusion map as we like by requiring that $J'$ is sufficiently close to $J$.

 In the case when
$q=n$ the result is essentially covered by Hamilton's work in \cite{kn:H}. Alternatively, still in the case when $q=n$, the result is essentially teh same as that proved by Tian in \cite{kn:T}, using an approach through the rigidity of the CR structure on $S^{2q-1}/\Gamma$. However the case $n>q$ seems to have essential new features, since the domain $V(a)$ does not then have a smooth boundary.

The \lq\lq rigidity'' expressed by Theorem 3 is related, in a more algebraic context, to the rigidity of quotient singularities proved by Schlessinger \cite{kn:S}. There is a notable distinction between the case $q\geq 3$ covered by the Theorem and the case $q=2$. In the latter case rigidity certainly does not hold. For example when $\Gamma=\pm 1$ the singularity $\bC^{2}/\Gamma$ is an affine quadric cone $Q(z)=0$ which can be deformed into a nonsingular quadric $Q(z)=\epsilon$. Correspondingly, the complex structure on the quotient of the annulus in $\bC^{2}$ is not rigid. The distinction between the cases
$q=2$ and $q\geq 3$ appears through the vanishing of $1$-dimensional sheaf cohomology in the latter case but not in the former. In fact it is very well known that products $\bC^{n-2}\times \bC^{2}/\Gamma$ can appear as Gromov-Hausdorff limits of blow-up sequences, under our hypotheses: the simplest example being when $n=2$. See also the further discussion in Section 5.

\section{Theorem 2 implies Theorem 1}

Suppose $(M_{i},g_{i})$ is a sequence of manifolds satisfying our standing conditions (with fixed $n,D$) and that $B(x_{i}, r_{i})$ are balls in $M_{i}$, so $r_{i}\leq D$. Suppose that the normalised energies $E(x_{i}, r_{i})$ tend to zero as $i\rightarrow \infty$.  Then the pointed manifolds $(M_{i}, x_{i})$ with the rescaled metrics $r_{i}^{-2}g_{i}$ have a Gromov-Hausdorff convergent subsequence, which we may as well suppose is the full sequence. Let $x_{\infty}, M_{\infty}$ be the based limit.   We claim that there are no singular points in the interior ball $B(x_{\infty},1/2)$.
This statement implies Theorem 1. For suppose Theorem 1 is false, so there is a sequence of balls $B(x_{i},r_{i})$ as above but point $y_{i}\in B(x_{i}, r_{i}/2)$ with $ \vert \Riem\vert(y_{i})=K_{i} r_{i}^{-2}$ with $K_{i}\rightarrow \infty$. We get a contradiction to the fact that the rescaled metrics converge in $C^{\infty}$ on the regular part of the limit $M_{\infty}$.

To prove the claim above we again argue by contradiction.  Notice that the the metric on the regular part of the unit ball $B(x_{\infty}, 1)$ in $M_{\infty}$ is flat. Suppose that, contrary to the claim,  $y\in B(x_{\infty},1/2)$ is a singular point. A tangent cone to $M_{\infty}$ at $y$ has the form $\bC^{n-q_{0}}\times C(Y_{0})$ for a length space $Y_{0}$.
If $Y_{0}$ is itself singular we take a tangent cone to $\bC^{n-q_{0}}\times C(Y)$ at a singular point and by the general Cheeger, Colding, Tian theory this must have the form $\bC^{n-q_{1}}\times C(Y_{1})$ for some $q_{1}<q_{0}$.
After at most $n$ steps  we arrive at an iterated tangent cone of the form $\bC^{n-q}\times C(Y)$ with $Y$ smooth. Since the metric on the regular part is flat we have $C(Y)= \bC^{q}/\Gamma$ where $\Gamma$ acts freely on the sphere. Passing to subsequences we can find a sequence of points $x'_{i}$ in $M_{i}$ and rescalings $l_{i}\rightarrow \infty$ such that $(M_{i}, x'_{i}, l_{i}^{2} g_{i}) $ have based limit $\bC^{n-q}\times \bC^{q}/\Gamma$.
Thus we can deduce from Theorem 2 that $q\leq 2$ and in fact, by the result of Cheeger \cite{kn:C}, the only possibility is $q=2$. 

To finish the proof, dealing with the singularities of complex codimension 2, we invoke the result of Cheeger,Colding and Tian from \cite{kn:CCT} which was also crucial in our previous paper \cite{kn:CD}.  Given any $\alpha>0$ we can find a ball $B(z,s)\subset M_{\infty}$ such that the Gromov-Hausdorff distance from $B(z,s)$ to the ball of radius $s$ in the model $\bC^{q}\times \bC^{q}/\Gamma$ is less than $\alpha s/2$. Now fix $s$ and choose $i$ so large that for a suitable choice of $x''_{i}\in M_{i}$ the Gromov-Hausdorff distance from $B(x''_{i}, s)$ to $B(z,s)$ is also less than $\alpha s/2$. Then Theorem 8.1 in \cite{kn:CCT} tells us that, for a suitable choice of $\alpha$, we have a fixed $\eta>0$ such that
$$   \int_{B(x'',s)} \vert \Riem \vert^{2} \geq \eta s^{n-2}. $$
This contradicts our hypothesis that $E(x_{i}, r_{i})\rightarrow 0$.

\section{Theorem 3 implies Theorem 2}

Consider a compact differentiable submanifold $\Sigma$, of dimension $2q-1$, in a K\"ahler manifold $M,\omega$. Suppose that $H^{1}(\Sigma, \bR)= H^{2}(\Sigma,\bR)=0$.
Thus we can write $\omega\vert_{\Sigma}= d\theta$ for a $1$-form $\theta$ on $\Sigma$ and the integral
\begin{equation}   I(\Sigma,\omega)= \int_{\Sigma} \omega^{q-2}\wedge \theta, \end{equation}
does not depend on the choice of $\theta$. 

Let $S\subset V_{\Gamma}(a)$ be the quotient of the unit sphere in $\{0\}\times \bC^{q}\subset \bC^{n-q}\times \bC^{q}$.

\begin{prop}
Let $M,\omega$ be a compact K\"ahler manifold of complex dimension $n$ such that $[\omega/2\pi R]$ is an integral class.
Suppose there is a holomorphic embedding $\iota:V_{\Gamma}(a)\rightarrow M$, for some $a>1$ and let $\Sigma=\iota(S)$. Then $I(\Sigma,\omega)\geq (2\pi)^{q} R^{q-1}$.  
\end{prop}

This is essentially standard complex geometry. By scaling, there is no loss in taking $R=1$. Let $L\rightarrow M$ be a holomorphic line bundle with curvature form $-i\omega$ and choose a power $k>0$ so that the sections of $L^{k}$ give an embedding of $M$ in $\bC\bP^{N}$. For $a>1$ let $W(a)$ be the annulus $\{w: a^{-1}<\vert w \vert<a\}$ in $\bC^{q}$ and $W_{\Gamma}(a)$ be the quotient by the free action of $\Gamma\subset U(q)$.
Let $\pi:W(a)\rightarrow W_{\Gamma}(a)$ be the quotient map. The pull-back
$\pi^{*}\iota^{*}(L^{k})$ is a holomorphic line bundle over $W(a)$. Since
$H^{1}(W(a); \cO)$ vanishes this line bundle is trivial. The line bundle
$\iota^{*}(L^{k})$ is determined by a character of $\Gamma$ so, increasing $k$ if necessary, we may suppose this is also trivial. Fix a trivialising section $\sigma$ of $\iota^{*}(L^{k})$. Thus the composite $W_{\Gamma}(a)\rightarrow
M\rightarrow \bC\bP^{N}$ is given by sections $s_{i}= f_{i}\sigma$ for $i=0,\dots N$, where $f_{i}$ are holomorphic functions on $W_{\Gamma}(a)$. 

Let $\Delta$ be the ball $\{ w: \vert w \vert <a^{-1}\}$ and  $\Delta_{\Gamma}$ be the quotient by $\Gamma$. Lifting the $f_{i}$ to $W(a)$ and applying Hartogs' theorem we see that they extend to holomorphic functions on $\Delta_{\Gamma}$
(i.e. $\Gamma$-invariant functions on $\Delta$). Thus $\iota$ extends to a holomorphic map $\iota^{+}: \Delta_{\Gamma}\setminus T\rightarrow M\subset \bC\bP^{N} $, where $T$ is a discrete subset defined by the common zeros of the extended functions $f_{i}$. Let $Z$ be  the graph of
$\iota^{+}$ and $\oZ$ be the closure of $Z$ in $\Delta_{\Gamma}\times \bC\bP^{N}$.  Thus $Z$ is an analytic variety and we have  holomorphic maps $p:\oZ\rightarrow \Delta_{\Gamma}$ and $j:\oZ\rightarrow M $. Writing $\sigma=  f_{i}^{-1} s_{i}$ we see that $\sigma$ defines a meromorphic section of $j^{*}(L^{k})$ with no zeros but with possible poles along a divisor supported in $\oZ\setminus Z$, corresponding to the points of $T$. 

We have to see that $T$ is nonempty, so that $\sigma$ does indeed have some poles. If $T$ is empty then $p:\oZ\rightarrow \Delta_{\Gamma}$ is a holomorphic equivalence so $j$ can be viewed as a  holomorphic map from $\Delta_{\Gamma}$ to $M$.  We extend the argument in the obvious way  to construct a  holomorphic map $J$ from $B\times \Delta_{\Gamma}$ to $M$, where $B$ is a ball in $\bC^{n-q}$, with $J$ equal to the the embedding $\iota$  on $B\times  V_{\Gamma}(a)$.
But it is clear that this is impossible if $M$ is smooth, as we suppose. 

Now we regard
$c_{1}(j^{*}L^{k})$ as a compactly supported cohomology class on $\oZ$, using the trivialisation $\sigma$ over the boundary. Since $\sigma$ has poles  we have
\begin{equation} \omega^{q-1} \wedge c_{1}(j^{*}L^{k})<0 . \end{equation}
The integrality of the Chern class then implies that
\begin{equation} \omega^{q-1}\wedge c_{1}(j^{*}L^{k})\leq -(2\pi)^{q-1} k. \end{equation}
Let $\theta$ be the $1$-form $k^{-1}\frac{i}{2}(\db-\partial)\log \vert \sigma\vert^{2}$ on $Z$. Then, regarding $\theta$ as a current on $\oZ$ we have an equation of currents
$$   d\theta= \omega + 2\pi k^{-1} E, $$
where $E$ is the current of the divisor representing $c_{1}(j^{*}L^{k})$, as a compactly supported cohomology class. If $\oZ_{0}\subset \oZ$ is the region interior to $\Sigma$, in the obvious sense, then by Stokes' Theorem,
$$ I(\Sigma,\omega)=  \int_{\Sigma} \theta \wedge \omega^{q-2}= -k^{-1}\omega^{q-1}\wedge c_{1}(j^{*}L^{k})+\int_{\oZ_{0}} \omega^{q} \geq (2\pi)^{q}. $$

\

\

It is now easy  to deduce Theorem 2 from Theorem 3 and the Proposition above.
Suppose that $x_{i}, M_{i}, l_{i}^{2} g_{i}$ is a sequence as considered in Theorem 2, with based Gromov-Hausdorff limit the length space $\bC^{n-q}\times \bC^{q}/\Gamma$.
Recall that we have the non-collapsing condition (2) and a two sided bound on the Ricci curvature of the $M_{i}, g_{i}$ by (4).
 By standard theory (using results of Anderson \cite{kn:MA}), this means that the metrics converge in $C^{\infty}$ on the smooth part of $\bC^{n-q}\times \bC^{q}/\Gamma$. Regarding the complex structures as covariant constant tensors we see that we can also suppose these converge. This means  that if we fix any $a>1$ we can find embeddings $\chi_{i}: V_{\Gamma}(a)\rightarrow M_{i}$ such that the pull-backs of the metrics and complex structures by $\chi_{i}$ converge to the standard structures on $V_{\Gamma}(a)$.  Applying Theorem 3 we see that we can suppose the $\chi_{i}$ are holomorphic embeddings.
So we are in  the situation considered in Proposition 1 with submanifolds $\Sigma_{i}$.  Applying Proposition 1 we see that $I(\Sigma_{i}, l_{i}^{2} \omega_{i})\geq (2\pi)^{q} l_{i}^{2q-2}\rightarrow \infty$. But this is a contradiction since $I(\Sigma_{i},l_{i}^{2} \omega_{i})$ is determined by the restriction of $l_{i}^{2}\omega_{i} $ to $\Sigma_{i}$, which converges to the standard model as $i\rightarrow \infty$.

\section{Proof of Theorem 3: Complex rigidity}

To simplify notation we will prove the result for some particular pair
$a_{1}, b_{1}$ but it will be clear that the argument can be adjusted to any pair. We will work with the domains $V(a)$; equivariance under the action of $\Gamma$ will allow us to deduce the result for the quotient spaces. We will consider various values of the parameter $a$, but all lying in some fixed interval, say $2\leq a\leq 4$.   

 According to \cite{kn:CS}  any (0,1) form $\sigma$ on $V(a)$ can be expressed
 as
$$  \sigma=\db K(\sigma) + K \db \sigma $$
where $K(\sigma)$ is orthogonal to the $L^{2}$ holomorphic functions and
$K(\db \sigma)$ is orthogonal to the image of $\db$. The   \lq\lq Kohn operator''
$K$ is bounded on $L^{2}$ (see \cite{kn:CS}, Section 2.3). It is easy to check  that this bound can be taken  independent of $a$. Notice that it as at this stage that the vanishing of $H^{1}(V(a); \cO)$ is fed into the proof.

Fix $\alpha \in (0,1)$ and let $\Vert\ \Vert_{k,\alpha,a}$ denote the $C^{k,\alpha}$ norm over $V(a)$. 
Now consider a different parameter $a^{*}<a$, so $V(a^{*})\subset V(a)$.   
\begin{prop} There are fixed $C,p$ such that
\begin{itemize}\item If $g$ is a function on $V(a)$ we have
$$   \Vert g \Vert_{2,\alpha,a^{*}} \leq C (a-a^{*})^{-p} \left( \Vert \db g\Vert_{1,\alpha,a}
+ \Vert g \Vert_{L^{2}(V(a)}\right). $$
\item If $\tau$ is a $(0,1)$-form on $V(a)$ with $\db^{*}\tau=0$ then
 $$   \Vert \tau \Vert_{1,\alpha,a^{*} } \leq C (a-a^{*})^{-p} \left( \Vert \db \tau \Vert_{0,\alpha,a}
+ \Vert \tau \Vert_{L^{2}(V(a)}\right). $$
\end{itemize}
\end{prop}

To see this, we can cover $V(a^{*})$ by balls of radius $(a-a^{*})/10$ say, such that for each ball the twice-sized ball with same centre is contained in $V(a)$. On a unit sized ball we have a standard elliptic estimate for functions
$$    \Vert g \Vert_{C^{1,\alpha}(B/2)}\leq {\rm const.} \left( \Vert \db g \Vert_{C^{0,\alpha}(B)} + \Vert g \Vert_{L^{2}(B)}\right) $$
and similarly for $(0,1)$ forms. Now the result follows by scaling. 

\

For the rest of the proof we will use the standard convention that $C,p$ are constants which may change from line to line. 

\

We are now ready to begin our main construction.  
 Consider a deformed complex structure on $V(a)$, defined by a tensor $\mu=\sum \mu_{ij} d\overline{z}_{j} \otimes \frac{\partial}{\partial z_{i}}$, smooth up to the boundary. Thus we have a deformed $\db$-operator $\db_{\mu}= \db +\mu \partial $. Let $f$ be a holomorphic function on $V(a)$, for the standard complex structure. Thus
$ \db_{\mu} f= \beta$ where $\beta=\mu \partial f$. (In our application $f$ will be one of the co-ordinate functions on $\bC^{n}$). Then $\beta= \db K(\beta) + K(\db \beta)$. Write $g=K(\beta)$ and $\beta'=K(\db \beta)$. Thus $\db g= \beta-\beta'$ while
$$   \db \beta'= \db \beta  \  \  \ \ \  \  \db^{*}\beta'=0. $$
The integrability of the deformed complex structure gives
$$   \db \beta = \db_{\mu}(\beta) - \mu \partial \beta = \db_{\mu}^{2} f - \mu\partial \beta = -\mu \partial \beta. $$
Applying the second item in the proposition above, and the $L^{2}$-boundedness of the Kohn operator, we get
\begin{equation}   \Vert \beta' \Vert_{1,\alpha,a^{*}} \leq C (a-a^{*})^{-p} \Vert \mu \partial \beta \Vert_{0,\alpha,a}. \end{equation}  
Applying the second item  to an intermediate region and then the first item we obtain

\begin{equation}  \Vert g \Vert_{2,\alpha,a^{*}}\leq  C (a-a^{*})^{-p}\left( \Vert \beta\Vert_{1,\alpha,a} +
 \Vert \mu \partial \beta \Vert_{0,\alpha,a}\right) \end{equation}
Now write $f'=f-g$. We have $\db_{\mu} f'= \beta' - \mu \partial g$
so
\begin{equation}  \Vert \db_{\mu} f' \Vert_{1,\alpha,a^{*}}\leq C (a-a^{*})^{-p} \left( \Vert \mu\Vert_{1,\alpha,a} \Vert \beta \Vert_{1,\alpha,a} + \Vert  \mu \Vert_{1,\alpha,a}^{2} \Vert \beta \Vert_{1,\alpha,a}\right).  \end{equation}

We make this construction starting with the $n$ co-ordinate functions
$f_{i}= z_{i}$, and getting new functions $f'_{i}$. Then $\beta_{i}=\mu \partial  f_{i}$ are just the components of $\mu$:
 \begin{equation}  \beta_{i}= \sum_{j} \mu_{ij} d\overline{z}_{j}. \end{equation}
  Let $\uf':V(a)\rightarrow \bC^{n}$ be the map with components $f'_{i}$. Suppose that the restriction of $\uf'$ is a diffeomorphism from $V(a^{*})$ to its image in $\bC^{n}$ and that the image contains a domain $V(a')$ where $a'$ is slightly less than $a^{*}$. Let $F:V(a')\rightarrow V(a^{*})$ be the inverse diffeomorphism. We  transport the complex structure defined by $\mu$ to $V(a')$, using the map $F$. Write
$$   \db_{\mu} f'_{i}= \sum \tau_{ij} d\overline{z}_{j}\ \ \ \ \ \  \ \ \db_{\mu} \overline{f'_{i}}= \sum D_{ij} d\overline{z}_{j}, $$
and suppose that the matrix $(D_{ij})$ is invertible at each point. A straightforward calculation shows that the \lq\lq new'' complex structure on $V(a')$ is defined by a tensor $\mu'$ which is given in matrix notation by
\begin{equation}  \mu'(z)= (D^{-1} \tau)(F(z)).\end{equation}

 The upshot is that, provided the various conditions above are met we get a complex structure defined by $\mu'$ on $V(a')$, given  by the formula (12), and a diffeomorphism
$F:V(a')\rightarrow V(a)$ which intertwines $\mu'$ and $\mu$.

\

We want to iterate this procedure, provided always that the initial deformation is sufficiently small.   We start by fixing a decreasing sequence of domains. Let $a_{1}=4$ and for integers $r\geq 2$ set 
$$  a_{r}= 4- \sum_{i=2}^{r} \frac{1}{i^{2}}, $$
which means that $a_{r}\geq  3$ for all $r$, so we take $b_{1} =3$.  Let $a^{*}_{r}= \frac{1}{2}(a_{r}+ a_{1})$, so $a_{r}-a^{*}_{r}= 1/2(r+1)^{2}$. Suppose we start with a $\mu_{1}$ on $V(a_{1})$ and that at stage $r$ we have constructed $\mu_{r}$ on $V(a_{r})$with a diffeomorphism $ {\cal F}_{r}: V(a_{r})\rightarrow V(a_{1})$ which intertwines
$\mu_{r}$ and $\mu_{1}$. Then, provided the various conditions above are met we perform the construction above to get $\mu_{r+1}$ on $V(a_{r+1})$ and a diffeomorphism  $F_{r+1}: V(a_{r+1})\rightarrow V(a_{r})$ so we can continue the inductive construction with ${\cal F}_{r+1}= {\cal F}_{r} \circ F_{r+1}$.

We need to show that, if $\Vert \mu_{1}\Vert_{1,\alpha,a_{1}}$ is sufficiently small then 
\begin{itemize}
\item The construction can proceed at each stage. 
\item The restriction of the $\mu_{r}$ to the fixed interior domain $V(b_{1})$ tend to zero in $\Vert \ \Vert_{1,\alpha, b_{1}}$.
\item The restrictions of the diffeomorphisms ${\cal F}_{r}$ to $V(b_{1})$ converge
in $C^{2,\alpha}$ to a diffeomorphism ${\cal F}:V(b_{1})\rightarrow V(a_{1})$, which can be made as close as we please to the inclusion map by assuming $\mu_{1}$ sufficiently small. 
\end{itemize}

If we establish these facts then we prove Theorem 3 as follows. Given a deformed complex structure on $V_{\Gamma}(a_{1})$ we lift it to a $\Gamma$-invariant structure on $V(a_{1} )$. It is clear that the $\mu_{r}$ we construct at each stage are $\Gamma$-invariant and the diffeomorphisms are $\Gamma$-equivariant.
Then ${\cal F}$ induces the desired   diffeomeorphism from $V_{\Gamma}(b_{1})$ to $V_{\Gamma}(a_{1})$.

Suppose we have constructed $\mu_{r}, {\cal F}_{r}$. Let $\uf_{r}: V(a_{r})\rightarrow
\bC^{n}$ be the map defined as above. Then the conditions for proceeding to the next stage will all be met if $\uf_{r}$ is sufficiently close to the identity in $C^{2,\alpha}$. By (9) this will be the case if $\mu_{r}$ is sufficiently small in $C^{1,\alpha}$. More precisely, if we write $\eta_{r}= \Vert \mu_{r} \Vert_{1,\alpha,a_{r}}$ then  we can proceed to the next stage if
\begin{equation} \eta_{r}\leq \epsilon r^{-p_{0}}, \end{equation}
for some suitable fixed $\epsilon,p_{0}$. Now we can estimate $\mu_{r+1}$ using
(10), (12) and the behaviour of H\"older norms under compositions and products.
We get \begin{equation} \eta_{r+1} \leq  r^{p} \eta_{r}^{2}. \end{equation} It is neater to express this as
\begin{equation} \eta_{r+1} \leq C_{1} \frac{r^{2p_{1}}}{(r+1)^{p_{1}}} \eta_{r}^{2}. \end{equation}
for some fixed $C_{1}, p_{1}$.
For then if we write $w_{r}= C_{1} r^{p_{1}} \eta_{r}$ we simply have
$w_{r+1}\leq w_{r}^{2}$. Choose $k>0$ so that
$$  \exp(-k2^{s-1})< C_{1} s^{p_{1}} \epsilon s^{-p_{0}}, $$
for all $s\geq 1$. Then if $w_{1}\leq e^{-k}$, that is to say if $\eta_{1}$ is sufficiently small,  it follows by induction that  the condition (13) is met at each stage and $w_{r} \leq  \exp(-k2^{r-1})$. Thus the iteration can proceed for all $r$ and it is clear that the other conditions itemised above are met, because of the very rapid decay of the $\eta_{r}$.

\section{Discussion}
\begin{enumerate}
\item In this paper we have concentrated on  proving what we need for our main result. However it seems  likely that the arguments in the proof of Theorem 2 can be extended to obtain a precise description of the complex structure for a K\"ahler-Einstein manifold close to a singular limit $\bC^{n-q}\times \bC^{q}/\Gamma$,  when we drop the integrality condition on the K\"ahler class. As Tian has suggested, one expects the complex structure in such a case to be a crepant resolution of the quotient singularity (at least when $\Gamma\subset SU(q)$), and one expects the metric to be modelled on Joyce's ALE metric. This seems significant because,  in the general Cheeger, Colding, Tian theory, rather little is known about the structure of Einstein metrics close to a singular limit.

\item Our rigidity result, Theorem 3, can clearly be generalised to other settings. (For example, we could consider any domain $V$ in $\bC^{n}$ such which has a suitable exhaustion by subsets $V_{a}$ such that 
$H^{1}(V_{a};\cO)=0$ and on which $\db$ has closed image.)  It  fits into a long line of similar statements, beginning with the Newlander-Nirenberg Theorem on the integrability of almost-complex structures. Our proof has some relation to the proof by Kohn \cite{kn:K} of this theorem, and also with the results of Hamilton \cite{kn:H} for more general domains. The solution of the \lq\lq $\db$-problem'' is an essential ingredient in all these results. However there is a notable difference in our case. In \cite{kn:K}, \cite{kn:H} it is first established that the $\db$-problem has a  solution, obeying suitable uniform estimates,  for all small deformations of the complex structure. In our  situation we do not have such a statement: we only know that the problem can be solved for the unperturbed solution using (essentially) the K\"unneth formula. This is the reason why we have to introduce the \lq\lq shrinking domains'' in the problem.

Note also that since we  allow ourselves to shrink the domain  we do not really need the full force of the result of Chakrabarti and Shaw, so the method may extend still further

\end{enumerate}


\end{document}